\documentclass[12pt]{article}

\usepackage{graphicx,amsmath,amssymb,url,enumerate,mathrsfs,epsfig,color}
\usepackage{tikz}

\usepackage{pgfplots}

\usetikzlibrary{calc, patterns,arrows, shapes.geometric}
\pgfplotsset{compat=1.8}

\def\QED{\hskip0.1em\hfill\null\ \null\nobreak\hfill
\kern3pt\lower1.8pt\vbox{\hrule\hbox
{\vrule\kern1pt\vbox{\kern1.7pt \hbox{$\scriptstyle
QED$}\kern0.2pt}\kern1pt\vrule}\hrule}}
\title{Material groupoids and algebroids}
\author{Marcelo Epstein\footnote{Department of Mechanical and Manufacturing Engineering, University of Calgary, Canada.} ~and
Manuel de Le\'on\footnote{Instituto de Ciencias Matem\'aticas
CSIC-UAM-UC3M-UCM, Spain}}
\date{}
\begin{document}
\maketitle
%%%%%%%%%%%%%%%%%%%%%%%%%%%%%%%%%%%%%%%%%%%%%%%%%%%%%%%%%%%%%%%%%%%%%%%%%%%%%%%%%%%%%%%%%%%%%%%%%%%%%%%%%%%%%
%%%%%%%%%%%%%%%%%%%%%%%%%%%%%%%%%%%%%%%%%%%%%%%%%%%%%%%%%%%%%%%%%%%%%%%%%%%%%%%%%%%%%%%%%%%%%%%%%%%%%%%%%%%%%

\bigskip
\begin{abstract}
Lie groupoids and their associated algebroids arise naturally in the study of the constitutive properties of continuous media. Thus, Continuum Mechanics and Differential Geometry illuminate each other in a mutual entanglement of theory and applications. Given any material property, such as the elastic energy or an index of refraction, affected by the state of deformation of the material body, one can automatically associate to it a groupoid. Under conditions of differentiability, this material groupoid is a Lie groupoid. Its associated Lie algebroid plays an important role in the determination of the existence of material defects, such as dislocations. This paper presents a rather intuitive treatment of these ideas.
\end{abstract}

{\bf Keywords}: Groupoids, algebroids, material uniformity, connections, material symmetries.

\section{Introduction}

Although not widely used, the notion of {\it material groupoid} is one of the most intuitive ideas in Continuum Mechanics. Given a material body $\mathcal B$ and considering some material property, we ask the following question: Do two points $X, Y \in {\mathcal B}$ share this particular material property? To answer this question we need to establish a {\it means of comparison}. This may consist, for instance, of a local diffeomorphism between neighborhoods of $X$ and $Y$. Moreover, we need to establish a means of evaluating whether or not the comparison has been successful. If it has, we draw and arrow from $X$ to $Y$. Otherwise, no arrow is drawn. Since `having the same property' is surely an {\it equivalence relation}, every point $X \in {\mathcal B}$ should be assigned a loop-shaped arrow. Moreover, for every arrow drawn from $X$ to $Y$ there should also be an arrow drawn from $Y$ to $X$. Finally, if there is an arrow from $X$ to $Y$ and another arrow from $Y$ to $Z$, there should also be an arrow from $X$ to $Z$. The set $\mathcal Z$ of all these arrows constitutes the {\it material groupoid associated with the chosen material property}. The body $\mathcal B$ is called the {\it base} of the groupoid. As a collection of arrows, the groupoid is undoubtedly a concrete geometric realization of the constitutive nature of a material body. For example, if every pair of points is joined by an arrow, the body is {\it uniform} with respect to the chosen property. In the mathematical terminology, we speak of a {\it transitive} groupoid. At the other extreme, if no two points are connected by an arrow (so that the only surviving arrows are the loops), we have a {\it totally intransitive groupoid}.

Clearly, these rough intuitive notions can be properly formalized; but, for now, let us assume that in some sense the arrows depend smoothly on the pairs of points. By this we mean that (thinking of a transitive groupoid, for definiteness) if $X$ is close to $X'$ and $Y$ is close to $Y'$ then the arrow from $X$ to $Y$ is also close to the arrow from $X'$ to $Y'$ and that this dependence of the arrow on the source-target pair is also differentiable, in some sense to be made precise later. We speak then of a {\it Lie groupoid}. We can look at the vicinity of a loop and calculate thereat the derivative of the arrow, that is, the linear part of the change in the arrow. In a manner completely and justifiably reminiscent of the passage from a Lie group to its Lie algebra (of `infinitesimal generators') we obtain a new geometric structure called the {\it Lie algebroid} associated with the Lie groupoid. The elucidation of the physical meaning of this construct is the main objective of this short paper.

In Section \ref{sec:groups} the concept of groupoid is briefly reviewed as a purely algebraic generalization of the notion of group. Section \ref{sec:alternative}, on the other hand, reviews an equivalent definition more amenable to geometric considerations. The definition of a Lie groupoid is provided within this context. In section \ref{sec:transitive} special attention is devoted to transitive groupoids and their relation to principal fibre bundles. In Section \ref{sec:algebroid} an attempt is made to motivate intuitively the notion of the Lie algebroid associated with a Lie groupoid. This idea is strengthened in Sections \ref{sec:coordinates} and \ref{sec:topological} by a discussion of coordinate expressions and different ways to organize the information carried by a groupoid. Connection in general as horizontal distributions in a fibre bundle and, in particular, connections in principal bundles are reviewed in Section \ref{sec:algebroidconnection} and correlated with vector bundle morphisms from the tangent bundle of the base of a transitive Lie groupoid into the associated Lie algebroid. The next three sections are devoted to show how all these mathematical notions appear naturally when specifying the constitutive response of a uniform material body, namely, a body that is made of the same material at all its points. Specifically, if the associated groupoid is a transitive Lie groupoid, the classical notion of material connection is elucidated in terms of the corresponding Lie algebroid.

\section{Groups and groupoids}
\label{sec:groups}

Recall that a {\it group} consists of a set $\mathcal G$ and two operations. The first, called {\it inversion}, is a unary operation
\begin{eqnarray} \label{eq11}
{\rm inv}: {\mathcal G} &\to& {\mathcal G} \nonumber \\
g &\mapsto& {\rm inv}(g)=g^{-1}.
\end{eqnarray}
The second operation is binary and is called {\it group multiplication} or {\it product}. It is usually indicated by apposition of the operands, that is,
\begin{eqnarray} \label{eq12}
{\rm prod}: {\mathcal G} \times {\mathcal G} &\to& {\mathcal G} \nonumber \\
(g,h) &\mapsto& {\rm prod}(g,h)=gh.
\end{eqnarray}
These two operations are assumed to enjoy the following properties:
\begin{enumerate}
\item \underline{Associativity}
\begin{equation} \label{12a}
(gh)k=g(hk)\;\;\forall g,h,k \in {\mathcal G}.
\end{equation}
\item \underline{Existence of a unit (or identity) element $e$}
\begin{equation} \label{eq13}
\exists e \in {\mathcal G}\;:\;eg=ge=g\;\;\forall g \in {\mathcal G}.
\end{equation}
\item \underline{Inversion consistency}
\begin{equation} \label{eq14}
gg^{-1}=g^{-1}g=e\;\;\forall g \in {\mathcal G}.
\end{equation}
\end{enumerate}
It is not difficult to prove that the group identity is unique and so is the inverse of each group element. Moreover,
\begin{equation} \label{eq15}
\left(g^{-1}\right)^{-1}=g\;\;\forall g \in {\mathcal G},
\end{equation}
and
\begin{equation} \label{eq16}
\left(gh\right)^{-1}=h^{-1}g^{-1}\;\;\forall g, h \in {\mathcal G}.
\end{equation}

A group $\mathcal G$ is {\it commutative} or {\it Abelian} if $gh=hg\;\;\forall g,h \in {\mathcal G}$. In a commutative group it is customary to refer to the group product as {\it sum} or {\it addition} and to denote it by $+$. Moreover, the unit element of a commutative group is usually denoted as $0$. The inverse of $g$ is represented as $-g$.

Just like a group, a groupoid is a set $\mathcal Z$ with two operations, also called {\it inverse} and {\it product}. There is, however, an essential difference in respect of the product operation, namely, the product in a groupoid is not defined for all pairs of elements of $\mathcal Z$. In other words, given two elements $x,y \in {\mathcal Z}$, the product $xy$ may or may not exist. The operations must satisfy
\begin{enumerate}
\item \underline{Associativity}
\begin{equation} \label{eq17}
\exists\; xy, yz \Leftrightarrow \exists\; x(yz) \Leftrightarrow \exists\; (xy)z \Rightarrow (xy)z=x(yz).
\end{equation}
\item \underline{Inversion}
\begin{equation} \label{eq18}
\forall\; x \in {\mathcal Z}\;\exists\; x^{-1},\;\; x^{-1}x,\;\; xx^{-1}\;\; \in {\mathcal Z}.
\end{equation}
\item \underline{Units or identities}
\begin{equation} \label{eq19}
\exists\; xy \;\Rightarrow \; xyy^{-1}=x,\;\;x^{-1}xy=y.
\end{equation}
\end{enumerate}
It can be shown that Equations (\ref{eq15}) and (\ref{eq16}) apply, mutatis mutandis, in the case of groupoids.

\section{An alternative definition}
\label{sec:alternative}

The definition of a groupoid just given is undoubtedly elegant and useful in that it puts the finger exactly on how a groupoid is deficient, with respect to a group, in terms of the non-existence of a legitimately binary product. In a groupoid, indeed, the product is not defined over the whole Cartesian product ${\mathcal Z}\times{\mathcal Z}$. As a result of this deficiency, we must forego the existence of a unique identity element and replace it by many such elements in their appropriate contexts, as expressed in Equation (\ref{eq19}). Although many examples of groupoids can be produced and the stated properties of the operations checked in each case, we are unlikely to emerge from these examples with a more intuitive picture of what a groupoid structure really means.

An alternative, equivalent, definition of a groupoid is more amenable to a geometric mental representation \cite{weinstein}. The key idea behind the new definition consists of introducing an auxiliary set $\mathcal M$, known as the {\it base} of the groupoid $\mathcal Z$, on which the groupoid effect is apparent. If we focus our attention on Equation (\ref{eq18}), which is part of the original definition, we observe that there is a distinguished subset of members of the groupoid, namely, those that can be expressed in the form $xx^{-1}$, where $x$ ranges over the whole set $\mathcal Z$.\footnote{
Not all expressions of the form $xx^{-1}$ are necessarily distinct. Indeed, let $z \in {\mathcal Z}$ be such that the products $zx$ and $zy$ are defined. Then, according to (\ref{eq19}), we have $zxx^{-1}=zyy^{-1}=z$. It follows (by multiplying to the left by $z^{-1}$) that $xx^{-1}=yy^{-1}$. Conversely, let $xx^{-1}=yy^{-1}$. Then, according to (\ref{eq17}), if $zx$ is defined, so is $zy$.} We define the base set as
\begin{equation} \label{eq22}
{\mathcal M}=\{xx^{-1}\; : \; x \in {\mathcal Z}\}.
\end{equation}
We now define two (surjective) {\it projection} maps
\begin{eqnarray} \label{eq23}
\alpha:{\mathcal Z} &\to& {\mathcal M} \nonumber \\
z &\mapsto& z^{-1}z,
\end{eqnarray}
and
\begin{eqnarray} \label{eq24}
\beta:{\mathcal Z} &\to& {\mathcal M} \nonumber \\
z &\mapsto& zz^{-1},
\end{eqnarray}
called, respectively, the {\it source map} and {\it target map} of the groupoid. Notice that

\begin{equation} \label{eq25}
\beta(z)=\alpha(z^{-1})\;\;\;\;\forall\;z\in {\mathcal Z}.
\end{equation}

Following the terminological lead, we can regard every element $z \in {\mathcal Z}$ as an {\it arrow} directed from its source to its target. It is now a rather straightforward matter to verify that the product $yx$ is well defined if, and only if, $\beta(x)=\alpha(y)$. Moreover, $\alpha(yx)=\alpha(x)$ and $\beta(yx)=\beta(y)$. In short, arrows are composed in a tip-to-tail fashion, as in a child's game. A groupoid can be seen, therefore, as a collection of arrows and counter-arrows between some pairs of points of a set $\mathcal M$, where some rules of the game must be satisfied.

More rigorously, we define a groupoid ${\mathcal Z}\rightrightarrows{\mathcal M}$ as a set of {\it arrows} (${\mathcal Z}$) over a set of {\it objects} (${\mathcal M}$), two ({\it projection}) surjective maps
\begin{equation}
\alpha: \mathcal Z \to \mathcal M\;\;\;{\rm
and}\;\;\;\beta: \mathcal Z \to \mathcal M
\end{equation}
called, respectively, the {\it source} and the {\it target} maps,
and an operation ({\it composition}) defined only for those
ordered pairs $(y,z) \in \mathcal Z \times \mathcal Z$ such that
\begin{equation}
\alpha(z)=\beta(y).
\end{equation}
This operation, indicated by the reverse apposition of
the operands, must be {\it associative}, that is, $(xy)z=x(yz)$, whenever the products are defined. Moreover, at each point $m \in {\mathcal M}$, there exists an {\it identity} $id_m$ such that $z\;id_m = z$ whenever $\alpha(z)=m$, and
$id_m\;z=z$ whenever $\beta(z)=m$. Finally,
for each $z \in \mathcal Z$ there
exists a (unique) {\it inverse} $z^{-1}$ such that
$zz^{-1}=id_{\beta(z)}$ and $z^{-1}z=id_{\alpha(z)}$.
The map
\begin{eqnarray} \label{eq:inclusion}
i: {\mathcal M} \to {\mathcal Z} \nonumber \\
m \mapsto id_m
\end{eqnarray}
is called the {\it object inclusion map}.

It is this shorter alternative definition that we will bear in mind from now on. A groupoid is {\it transitive} if every pair of points is connected by at least one arrow. A groupoid is a {\it Lie groupoid} if both $\mathcal Z$ and $\mathcal M$ are differentiable manifolds, the projections are submersions, the inclusion map is smooth and so are the operations of composition and inversion.

The set ${\mathcal Z}_m = \alpha^{-1} (m)$ is called the $\alpha$-{\it fibre} over $m \in {\mathcal M}$. It consists of all the arrows issuing from $m$. Analogously, the $\beta$-{\it fibre} over $m \in {\mathcal M}$ is the set ${\mathcal Z}^m = \beta^{-1} (m)$ consisting of all arrows arriving at $m$. The intersection ${\mathcal Z}_m^n={\mathcal Z}_m \cap {\mathcal Z}^n$ consists of all arrows starting at $m$ and ending at $n$. Each of the sets ${\mathcal G}_m={\mathcal Z}_m^m$ is clearly a group, called the {\it vertex group at} $m$. If ${\mathcal Z}_m^n \ne \emptyset$, the corresponding vertex groups are necessarily conjugate, namely,
\begin{equation}
{\mathcal G}_n= z\, {\mathcal G}_m \, z^{-1},
\end{equation}
where $z \in {\mathcal Z}_m^n$. Thus, in a transitive groupoid all the vertex groups are mutually conjugate. Any one of them may be regarded as the {\it structure group} of the transitive groupoid.

\section{Transitive groupoids and principal bundles}
\label{sec:transitive}

A useful analogy to bear in mind when thinking of the relation between transitive groupoids and principal bundles is the relation between affine and vector spaces. Given an affine space, if we arbitrarily adopt a particular fixed point as reference, the affine space becomes essentially equivalent to a vector space. In this sense, we may say that the affine space of departure can be regarded as an equivalence class of vector spaces. Vice versa, given a vector space, we can regard it as an affine space by associating to any pair of vectors their difference.

In this spirit, let ${\mathcal Z}\rightrightarrows{\mathcal M}$ be a transitive Lie groupoid and let $m \in {\mathcal M}$ be an arbitrary point of the base. The disjoint union
\begin{equation}
{\tilde {\mathcal Z}}_m = \bigcup_{x \in \mathcal M}\, {\mathcal Z}_m^x
\end{equation}
deserves further attention. Firstly, we observe that the vertex group ${\mathcal G}_m$ has a natural right action $\mathcal R$ on ${\tilde {\mathcal Z}}_m$ defined by
\begin{equation}
{\mathcal R}_g(z) = z g\;\;\;\;\;\forall g \in {\mathcal G}_m,\;\;\;\;\;\forall z \in {\tilde {\mathcal Z}}_m,
\end{equation}
where on the right-hand side we apply the groupoid composition, with the understanding that $g\in {\mathcal G}_m = {\mathcal Z}_m^m \subset {\mathcal Z}$ and $z \in {\tilde {\mathcal Z}}_m \subset  {\mathcal Z}$. This action is obviously {\it fibre preserving} in the sense that
\begin{equation}
z \in {\mathcal Z}_m^x \; \Rightarrow \;  z g \in {\mathcal Z}_m^x.
\end{equation}
Moreover, this right action is {\it free}. Recall that a group action is free if the condition ${\mathcal R}_g(z)=z\;\;{\rm for\;some\;} z$ implies that $g$ is the group identity. In our case, this property follows directly from the existence of inverses in the original groupoid. These conditions\footnote{Local triviality is also needed, a condition that some authors assume to be part of the Lie groupoid definition.} guarantee that the set ${\tilde {\mathcal Z}}_m$ is a principal bundle with structure group ${\mathcal G}_m$ and bundle projection $\beta$. Had we started from a different point $n \in {\mathcal M}$, we would have obtained a different principal bundle ${\tilde {\mathcal Z}}_n$, with structure group ${\mathcal G}_n$, conjugate to ${\mathcal G}_m$. Thus, just like in the case of an affine space and its associated vector spaces, the groupoid can be captured by any of the conjugate principal bundles obtained in the manner described above by an arbitrary choice of `origin'. The converse passage from a principal bundle to a groupoid is also possible \cite{language}.

\section{The Lie algebroid of a Lie groupoid}
\label{sec:algebroid}

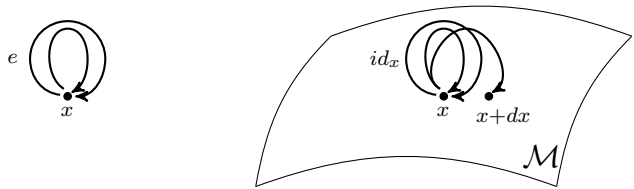
\begin{figure}[h!]
\begin{center}
\begin{tikzpicture} [scale=1.0]

\draw[fill](0,0) circle[radius=0.05cm];
\node[below]  at (0,0) {$_x$};
\draw [thick, stealth'-] (0.1,0) to [out=0, in=-90] (0.5,0.5) to [out=90,in=0] (0,1) to [out=180,in=90] (-0.5,0.5) to [out=-90, in=174] (-0.1,0.02);
\draw [thick, -stealth'] (-0.05,0.1) to [out=150, in=180]  (0,0.9) to [out=0,in=30]  (0.05,0.05);
\node[left] at (-0.5,0.5) {$_e$};
\begin{scope}[shift={(5,0)}]
\draw[fill](0,0) circle[radius=0.05cm];
\node[below]  at (0,0) {$_x$};
\draw[fill](0.6,0) circle[radius=0.05cm];
\node[below]  at (0.80,0) {$_{x+dx}$};
\draw [thick, stealth'-] (0.1,0) to [out=0, in=-90] (0.5,0.5) to [out=90,in=0] (0,1) to [out=180,in=90] (-0.5,0.5) to [out=-90, in=174] (-0.1,0.02);
\draw [thick, -stealth'] (-0.05,0.1) to [out=150, in=180]  (0.3,0.9) to [out=0,in=30]  (0.64,0.05);
\draw [thick, -stealth'] (-0.05,0.1) to [out=150, in=180]  (0,0.9) to [out=0,in=30]  (0.05,0.05);
\node[left] at (-0.4,0.5) {$_{id_x}$};

\end{scope}
\begin{scope}[shift={(2.5,-1.2)}, scale=2]
\draw (0,0) to [out=20, in=160] (2,0);
\draw (0,0) to [out=80, in=225] (0.5,1);
\draw (0.5,1) to [out=20, in=160] (2.5,1);
\draw (2,0) to [out=80, in=225] (2.5,1);

\node at (1.9,0.2) {$\mathcal M$};
\end{scope}

\end{tikzpicture}
\end{center}
\caption{Thinking in a Lie algebra (left) and in a Lie algebroid (right)}
\label{fig:algebroid}
\end{figure}

From the vantage point of the theory of groupoids a group $\mathcal G$ can be regarded as the particular case of a groupoid whose base is a singleton $\{x\}$. In this picture, therefore, the group is a collection of loop-shaped arrows starting and ending at $x$. Let us assume that $\mathcal G$ is a Lie group, namely, a group which is also a differentiable manifold where the operations of inversion and composition are smooth. These properties of a Lie group are so strong that, as demonstrated by Lie himself, the whole group can be understood by studying the behavior of the group elements in a small neighborhood of the unit element $e$. This observation can be formalized under the rubric of the {\it Lie algebra} associated with the Lie group. The corresponding construction for a Lie groupoid is the {\it Lie algebroid}. To anticipate intuitively the similarities and differences between the the Lie algebra of a Lie group and the Lie algebroid of a Lie groupoid, consider, as shown in Figure \ref{fig:algebroid}, that:
\begin{enumerate}
\item Whereas a Lie group has a single unit, a Lie groupoid has a unit associated with each and every point of the base manifold. In other words, in a Lie groupoid we have a field of units defined over the base manifold.
\item Therefore, while an element in a first-order vicinity of the group unit is a vector, an element in a first-order vicinity of the field of units is a {\it vector field}.
\item When exploring the vicinity of the Lie group unit, we naturally encounter group elements. Geometrically, we may say that these are all loop-shaped arrows obtained by keeping the end points coincident and fixed (since we have just one base point $x$) and deforming the loop. On the other hand, in a Lie groupoid, we can both deform the unit loop at $x$ (thus moving within the local group) and we can also open it by displacing, say, the tip of the arrow to a nearby point $x+dx$. This last feature constitutes the crucial difference between the Lie algebra and the Lie algebroid.
\end{enumerate}

From these considerations one can see that in the Lie algebroid $\mathfrak z$ of a Lie groupoid $\mathcal Z$ there will be elements within the local Lie algebras (`vertical' elements), namely, those obtained by moving strictly within arrows whose source and target are equal to some $x \in {\mathcal M}$. On the other hand, all other elements will be associated to a specific `opening' of the loop in some direction $dx$ lying in the tangent space of the base manifold $\mathcal M$. Clearly, then, the Lie algebroid (which is, obviously, a vector bundle over $\mathcal M$), projects naturally over the tangent bundle $T\mathcal M$. In other words, there exists a natural map $\rho: {\mathfrak z} \to {T\mathcal M}$ called the {\it anchor map} of the Lie algebroid. When the Lie groupoid is transitive, this map is surjective. Moreover, if the typical symmetry group of the transitive groupoid is trivial, then the anchor map is also one-to-one. Notice that what we called above vertical elements of the algebroid are precisely those that are anchored to the zero vectors of $T\mathcal M$.

Somewhat more technically, the Lie algebroid $\frak z$ of the transitive Lie groupoid ${\mathcal Z} \rightrightarrows {\mathcal M}$ is a vector bundle ${\frak z}=A{\mathcal Z}$ over $\mathcal M$ consisting of the disjoint union of the tangent spaces to the $\alpha$ fibres at the respective identities, namely,
\begin{equation}
A{\mathcal Z} = \bigcup_{x \in {\mathcal M}}  T_{id_x} {\mathcal Z}_x .
\end{equation}
Each of these tangent spaces consists of the vectors intuitively described on the right side of Figure \ref{fig:algebroid}. Properly speaking, the elements of the algebroid consist of sections of this vector bundle. The Lie bracket of the algebroid is the restriction of the ordinary bracket of sections of $T{\mathcal Z}$.

\section{Coordinate expressions}
\label{sec:coordinates}

\begin{figure}[h!]
\begin{center}
\begin{tikzpicture} [scale=1.0]

% spiders
\foreach \y in {0,1}
\foreach \x in {0,1,2}
{
\begin{scope}[shift={(4+2*\x+0.75*\y,1.5*\y+0.4*\x*(2.-\x))},scale=0.8]
%\draw [thick, stealth'-] (0.1,0) to [out=0, in=-90] (0.5,0.5) to [out=90,in=0] (0,1) to [out=180,in=90] (-0.5,0.5) to [out=-90, in=174] (-0.1,0.02);
\draw [thick, -stealth'] (-0.05,0.1) to [out=40, in=140]    (1.0,0.05);
\draw [thick, -stealth'] (-0.05,0.1) to [out=140, in=40]    (-1.0,0.05);
\draw [thick, -stealth'] (-0.05,0.1) to [out=80, in=160]    (0.75,0.45);
\draw [thick, -stealth'] (-0.05,0.1) to [out=100, in=20]    (-0.75,0.45);
\draw [thick, -stealth'] (-0.05,0.1) to [out=20, in=120]    (0.5,-0.3);
\draw [thick, -stealth'] (-0.05,0.1) to [out=160, in=60]    (-0.7,-0.3);
\draw [thick, -stealth'] (-0.05,0.1) to [out=110, in=160]    (0.35,0.6);
\draw [thick, -stealth'] (-0.05,0.1) to [out=190, in=120]    (-0.1,-0.4);
\draw [thick, -stealth'] (-0.05,0.1) to [out=150, in=180]  (-0.05,0.8) to [out=0,in=30]  (-0.05,0.1);
\end{scope}
}
% base manifold
\begin{scope}[shift={(2.6,-1.0)}, scale=3]
\draw (0,0) to [out=20, in=160] (2,0);
\draw (0,0) to [out=80, in=225] (0.5,1);
\draw (0.5,1) to [out=20, in=160] (2.5,1);
\draw (2,0) to [out=80, in=225] (2.5,1);
\node at (1.9,0.14) {$\mathcal M$};
\end{scope}

\end{tikzpicture}
\end{center}
\caption{The $\alpha$-bundle ${\mathcal Z}_{\mathcal M}$ as a spider colony}
\label{fig:spiders}
\end{figure}
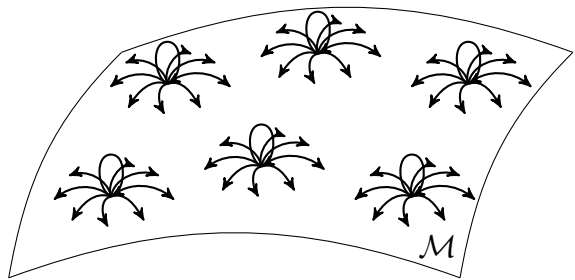

Although not strictly necessary, the manifestation of geometrical structures in terms of their expression in coordinate charts is often useful not only to perform actual calculations but also to reveal the close relation between apparently disparate constructs. Let us start by considering the disjoint union ${\mathcal Z}_{\mathcal M}$ of all the $\alpha$-fibres of a transitive Lie groupoid $\mathcal Z$, that is,
\begin{equation}
{\mathcal Z}_{\mathcal M}=\bigcup_{x \in \mathcal M} {\mathcal Z}_x.
\end{equation}
This set, which we call the $\alpha$-bundle, can be regarded as a fibre bundle over the base manifold $\mathcal M$ with projection $\alpha$. In terms of arrows, ${\mathcal Z}_{\mathcal M}$ looks like a spider colony, each fibre ${\mathcal Z}_x$ being a spider with legs issuing from $x$ and pointing to some point $y \in \mathcal M$, as shown schematically in Figure \ref{fig:spiders}. Notice that the total set of this fibre bundle is the same as the total set of the original transitive groupoid $\mathcal Z$. They both consist of the set of all arrows.

 Let a coordinate patch be given in $\mathcal M$ with coordinates $x^i,\;(i=1,...,\dim {\mathcal M})$. Let the structure group of $\mathcal Z$ be a Lie group $\mathcal G$ and let $\lambda^\xi\;(\xi=1,..., \dim{\mathcal G})$ be coordinates within this group such that $\lambda^\xi=0$ corresponds to the group identity. Then, an element of the bundle ${\mathcal Z}_{\mathcal M}$ (that is, a spider leg) is obtained by specifying a triple $\{(x^1,...),(\lambda^1,...),(y^1,...)\}$. Clearly, the end-point of the arrow (with coordinates $y^i$) may fall without the coordinate patch, unless the base manifold is trivial and has been covered with a single chart. This point is irrelevant for our construction, since we are going to be concerned with tangent vectors only. Indeed, an element in the tangent space to a fibre at the identity arrow, that is an element of $T_{id_x} {\mathcal Z}_x$, is a vector of the form
\begin{equation} \label{eq:coord1}
{\bf v}=a^\xi\;\left.\frac{\partial}{\partial \lambda^\xi}\right|_{\lambda=0,y=x}
 + \;\; b^i\;\left.\frac{\partial}{\partial x^i}\right|_{\lambda=0,y=x},
\end{equation}
 with some obvious notational abuse. Notice that the anchor map $\rho$ assigns to it the vector
 \begin{equation}
 \rho({\bf v})\;=\;b^i\;\left.\frac{\partial}{\partial x^i}\right|_{x}
 \end{equation}
 in $T_x{\mathcal M}$.

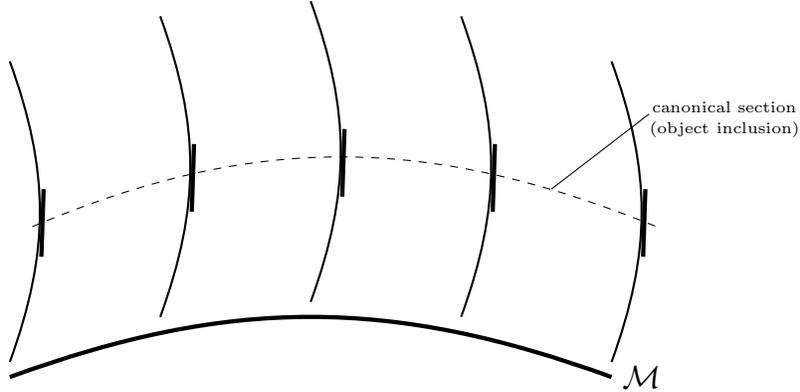
\begin{figure}[h!]
\begin{center}
\begin{tikzpicture} [scale=1.0]

% base manifold
\draw [ultra thick] (0,0) to [out=20, in=160]    (8,0);
% fibres
\foreach \x in {0,1,2,3,4}
{\begin{scope}[shift={(2*\x, 0.2-0.2*\x*(\x-4)}]
\draw [thick] (0,0) to [out=70, in=-70]    (0,4);
\draw[ultra thick] (0.42,1.4) to (0.45,2.3);
\end{scope}
}
%section
\draw [dashed] (0.3,2) to [out=23, in=158]    (8.6,2);
\node at (9.5,3.6){\tiny canonical section};
\node at (9.5,3.3){\tiny (object inclusion)};
\draw (8.5,3.5) -- (7.2,2.5);
\node at (8.4,0) {$\mathcal M$};
\end{tikzpicture}
\end{center}
\caption{The algebroid vector bundle (indicated with thicker lines)}
\label{fig:AZ}
\end{figure}

 Each of the fibres of ${\mathcal Z}_{\mathcal M}$ (that is, each spider) is , in fact, a manifold of dimension $\dim {\mathcal G} + \dim {\mathcal M}$. Moreover, this bundle has a canonical global section, which is essentially the object inclusion map (\ref{eq:inclusion}). The vector bundle $A{\mathcal Z}\rightarrow {\mathcal M}$ is schematically depicted in Figure \ref{fig:AZ}
  as a very particular sub-bundle of the tangent bundle
 $T{\mathcal Z}_{\mathcal M}$, by selecting out of each fibre (`spider') of ${\mathcal Z}_{\mathcal M}$ the tangent plane (to the fibre) at the identity.

\section{Topological considerations}
\label{sec:topological}

We have introduced several essentially equivalent topological (and differentiable) constructs. They differ mainly in the way they organize the information. Our purpose now is to suggest that there exist certain natural homeomorphims (and diffeomorphisms) between appropriately chosen neighbourhoods in a transitive Lie groupoid $\mathcal Z$, any one of its associated principal bundles ${\tilde{\mathcal Z}}_m$ and the spider-bundle
${\mathcal Z}_{\mathcal M}$. Recall that the elements of the principal bundle are arrows issuing from a single reference point $m \in {\mathcal M}$. These arrows are divided into fibres (equivalence classes) according to the value of the terminal point. In other words, two elements $z_1$ and $z_2$ are in the same fibre if, and only if, $\beta(z_1)=\beta(z_2)$. The structure group ${\mathcal G}_m$ acts on the right on the whole bundle. The elements in the spider bundle (which is not a principal bundle), on the other hand, comprise all the arrows of the original groupoid, but they are organized into fibres according to commonality of the $\alpha$-projection. The fibre over $m$ consists of the same arrows that are found in the totality of the principal bundle  ${\tilde{\mathcal Z}}_m$. We may say that the groupoid has been reorganized in such a way that each fibre contains exactly the arrows that make up one of the equivalent principal bundles! Each principal bundle has metamorphosed into a spider. Thus, it is not unreasonable to expect that there exists a natural continuous bijection between neighbourhoods in ${\tilde{\mathcal Z}}_m$ and respective neighbourhoods in the corresponding fibres of the spider bundle. The continuity is established within the bundle topology of ${\tilde{\mathcal Z}}_m$ and the subset topology induced in the spider fibres by the spider bundle topology. These natural homeomorphisms are, in fact, diffeomorphisms.

Let $x$ be a point in ${\mathcal M}$, and consider an open neighbourhood ${\mathcal U}_x \subset \alpha^{-1}(x)$ of $id_x \in \alpha^{-1}(x)$ within the fibre with the subset topology. Let
$z \in {\tilde{\mathcal Z}}_m$ be a point in an associated principal bundle (i.e. an arrow issuing from $m \in {\mathcal M}$) such that $\beta(z)=x$. Define the set
\begin{equation}
{\mathcal V}_z=\{v \in {\tilde{\mathcal Z}}_m\;|\; v = u z,\; u \in {\mathcal U}\}.
\end{equation}
A subtle point here is that the product $zu$ is defined only in the original groupoid $\mathcal Z$ and that, once this product is calculated, we obtain an arrow $v$ issuing from $m$ and, therefore, a legitimate element of the associated bundle ${\tilde{\mathcal Z}}_m$. The set ${\mathcal V}_z$ is an open set containing $z$ and, therefore, a legitimate bundle neighbourhood of $z$. The map ${\mathcal U}_x \to {\mathcal V}_z$ just constructed is also one-to-one and onto and, consequently, a bijection. Having been established without the use of coordinates, this bijection is canonical. Moreover, let ${\hat z} \in \beta^{-1}(x)$ be another element in the same fibre of ${\tilde{\mathcal Z}}_m$. Since the right action of ${\mathcal G}_m$ on ${\tilde{\mathcal Z}}_m$ is free (and transitive on fibres), there is a unique $g \in {\mathcal G}_m$ such that
\begin{equation}
{\hat z}= z g.
\end{equation}
Then
\begin{equation}
{\mathcal V}_{\hat z}=\{v \in {\tilde{\mathcal Z}}_m\;|\; v = u {\hat z}= u z g,\; u \in {\mathcal U}\}={\mathcal R}_g ({\mathcal V}_z).
\end{equation}
In other words, the canonical neighbourhoods on points on the same fibre of the associated principal bundle are obtained by right translation with an element of the structure group. Put differently, a collection of fibre-wise neighbourhoods of the identities in the $\alpha$-bundle of a given transitive topological groupoid $\mathcal Z$ uniquely (and canonically) determines an open cover of the associated principal groupoid ${\tilde{\mathcal Z}}_m$.

\section{Connections}
\label{sec:algebroidconnection}

 A {\it connection} in a fibre bundle $\pi: {\mathcal C} \to {\mathcal M}$ is a smooth {\it horizontal distribution} in $\mathcal C$. To make this definition more precise, for each $c \in {\mathcal C}$ let ${\mathcal V}_c$ denote the tangent space at $c$ to the fibre $\pi^{-1}\{\pi(c)\}$. A connection is a smooth assignment, at each point $c \in \mathcal C$, of a subspace ${\mathcal H}_c$ of the tangent space $T_c{\mathcal C}$ such that
 \begin{equation}
 T_c{\mathcal C}= {\mathcal H}_c \oplus {\mathcal V}_c,
 \end{equation}
 where $\oplus$ denotes the direct sum of vector spaces. It follows from this definition that the dimension of ${\mathcal H}_c$ is necessarily equal to $\dim {\mathcal M}$. All vectors $\bf w$ in ${\mathcal V}_c$ satisfy the condition
 \begin{equation}
 \pi_*({\bf w})={\bf 0}.
 \end{equation}
 For this reason, the subspaces ${\mathcal V}_c$ are called {\it vertical}. Correspondingly, the subspaces of the distribution are called {\it horizontal}.

The bundle projection (restricted to the horizontal subspace)
\begin{equation}
\pi_*: {\mathcal H}_c \to T_{\pi(c)} {\mathcal M}
\end{equation}
is an isomorphism of vector spaces. Its inverse
\begin{equation}
\Gamma: T_{\pi(c)} {\mathcal M} \to {\mathcal H}_c
\end{equation}
will be called the {\it Christoffel map} at $c$. Clearly, a horizontal distribution is equivalent to the specification of a smooth family of Christoffel maps. In a local coordinate system of the fibre bundle $\{(x^1, ...), (\lambda^1,...)\}$, the Christoffel maps are determined by (rectangular) matrices with entries $\Gamma^\xi_i = \Gamma^\xi_i (x^1, ...,\lambda^1,...)$. These entries are known as the {\it Christoffel symbols of the connection}. Specifically, for a vector ${\bf v} \in T_{\pi(c)}{\mathcal M}$ with components $v^i$, the image at $c$ by the Christoffel map is given by
\begin{equation}
\Gamma_c({\bf v})= v^i \frac{\partial}{\partial x^i} - \Gamma^\xi_i\, v^i \frac{\partial}{\partial \lambda^\xi}.
\end{equation}

A {\it principal bundle connection} is a connection in a principal bundle whose horizontal spaces are invariant under the structure group right action, namely,
\begin{equation}
\left({\mathcal R}_g\right)_*\,{\mathcal H}_c = {\mathcal H}_{{\mathcal R}_g c},\;\;\;\forall g \in {\mathcal G},\;\;c \in {\mathcal C}.
\end{equation}
This invariance condition has an important repercussion on the variation of the Christoffel maps along a fibre. Without the need of writing it explicitly, it is clear that for a principal bundle connection it is sufficient to know the Christoffel map at one point of a fibre in order to determine (by right translation) the Christoffel map at all other points of that fibre.

Suppose that a transitive Lie groupoid $\mathcal Z$ is at hand and that the corresponding Lie algebroid $\frak z= A{\mathcal Z}$ has been constructed. Suppose, moreover, that a vector bundle morphism
\begin{equation}
\gamma: T{\mathcal M} \to A{\mathcal Z}
\end{equation}
 has been given satisfying the condition
 \begin{equation}
 \rho \circ \gamma = id_{T{\mathcal M}}.
 \end{equation}
where $\rho$ is the anchor map of the Lie algebroid. In the coordinate system used in Equation (\ref{eq:coord1}), the vector bundle morphism $\gamma$ can be written as
\begin{equation}
u^i\;\mapsto\;\{\gamma^\xi_i u^i,u^i\},
\end{equation}
where $\gamma^\xi_i=\gamma^\xi_i(x^j)$ are functions defined on the base manifold. If we now consider the principal bundle ${\tilde{\mathcal Z}}_m$ defined in Section \ref{sec:transitive} and if we make use of the canonical diffeomorphisms introduced in Section \ref{sec:topological}, it is clear that we obtain a group-invariant horizontal distribution, that is, a legitimate principal bundle connection on ${\tilde{\mathcal Z}}_m$. Thus, there is a one-to-one correspondence between vector bundle morphisms from the tangent bundle of the base into the Lie algebroid and principal bundle connections in any one of the associated principal bundles ${\tilde{\mathcal Z}}_m$.

\section{Material groupoids}

Two points, $X$ and $Y$, of a simple material body $\mathcal B$ are {\it materially isomorphic} \cite{noll} if there exists a linear isomorphism $P_{XY}$ between the corresponding tangent spaces, that is,
\begin{equation} \label{eq100}
P_{XY}: T_X{\mathcal B} \to T_Y{\mathcal B},
\end{equation}
that brings their constitutive responses into coincidence. What is meant by `into coincidence' depends on the context. For example, in the case of a simple elastic body with constitutive equation
\begin{equation} \label{eq101}
{\bf T}={\bf T}({\bf F}, X),
\end{equation}
where $\bf T$ is the Cauchy stress and $\bf F$ is the deformation gradient, the condition may read
\begin{equation} \label{eq102}
{\bf T}({\bf F}, Y)={\bf T}({\bf FP}_{XY}, X).
\end{equation}

Every point $X$ is trivially materially isomorphic to itself via the identity transformation $\bf I$, but there may exist additional, non-trivial, material automorphisms. These are none other than the {\it material symmetries} of the constitutive law at $X$, which form the {\it material symmetry group} ${\mathcal G}_X$ at $X$. If $X$ and $Y$ are materially isomorphic, their symmetry groups are necessarily conjugate. Specifically, if $P_{XY}$ is a material isomorphism, then
\begin{equation} \label{eq103}
{\mathcal G}_Y={\bf P}_{XY}\;{\mathcal G}_X\;{\bf P}_{XY}^{-1}.
\end{equation}
Conversely, the degree of freedom in the choice of material isomorphisms between two materially isomorphic points is precisely the symmetry group, according to
\begin{equation} \label{eq104}
{\mathcal P}_{XY}={\bf P}_{XY}\;{\mathcal G}_X={\mathcal G}_Y\;{\bf P}_{XY}={\mathcal G}_Y\;{\bf P}_{XY}\;{\mathcal G}_X,
\end{equation}
where ${\mathcal P}_{XY}$ denotes the collection of all possible material isomorphisms between $X$ and $Y$.

Assigning an arrow to each material isomorphism between points of $\mathcal B$ we obtain the {\it material groupoid} $\mathcal Z$ of $\mathcal B$ relative to the given constitutive law \cite{jgp}. We have already painted this picture informally in the Introduction. The only difference now is that, because of the possible non-triviality of the symmetry groups, there may be many arrows between two points. The (never empty) collection of all the loop shaped arrows at a point $X$ is the geometrical representation of the symmetry group at $X$.

If the material groupoid is transitive, the body is said to be {\it materially uniform}. Since in a transitive groupoid all the local symmetry groups are conjugate, we can speak of a {\it typical symmetry group} $\mathcal G$ of the groupoid. In everyday language, a body is materially uniform if all its points are made of the same material. The acid test for this condition is that a small neighborhood of every point can be conceived as a perfect graft of some standard piece of material, with symmetry group $\mathcal G$, which we may call the {\it material archetype}. The material isomorphism represents what one has to do to this standard piece (rotation, deformation) before inserting it in place.\footnote{The notion of material isomorphism can be extended to encompass functionally graded materials \cite{fgm}.}

Since the body $\mathcal B$ is assumed to be a differentiable manifold, it is possible to consider the case in which the body is {\it smoothly materially uniform} in the sense that the material isomorphisms can be chosen as smoothly varying over each neighborhood. In that case, we speak of a {\it Lie groupoid}, whereby the base and the total set of arrows are differentiable manifolds and the projections are smooth surjective submersions and the groupoid operations are smooth in their arguments. Let $\mathcal B$ be a smoothly materially uniform body with Lie groupoid $\mathcal Z$. Adopting a reference configuration (and using a system of Cartesian coordinates $X^I$), the material isomorphisms boil down to matrices. Smoothness implies that in any given neighborhood these matrices can be chosen, within the degree of freedom afforded by the symmetry groups, so that their entries depend smoothly on position. If ${\bf P}(X)$ represents a material isomorphism from the archetype to $X \in {\mathcal B}$, the entries of the corresponding matrix are some (locally) smooth functions $P^I_\alpha(X^1, X^2, X^3)$. The Greek subscripts ($\alpha=1,2,3$) refer to components in a fixed Cartesian basis in the archetype.

\section{Material G-structures}

The material groupoid of a body $\mathcal B$ as defined in the previous section is always a subgroupoid of the {\it general Lie groupoid} $GL({\mathcal B})$ of $\mathcal B$ consisting of all the linear isomorphisms between the tangent spaces of all pairs of points of $\mathcal B$. According to the construction explained in Section \ref{sec:transitive}, any of the principal bundles associated with $GL({\mathcal B})$ is nothing but the {\it principal frame bundle} $F{\mathcal B}$, whose structure group is the general linear group $GL(3; {\mathbb R})$. Consequently, a principal bundle associated with a transitive Lie subgroupoid of $GL({\mathcal B})$ is a principal subbundle of $F{\mathcal B}$ with a structure group ${\mathcal G} \subset GL(3; {\mathbb R})$. In other words, a principal bundle associated with a transitive material Lie groupoid is precisely a {\it G- structure} \cite{sternberg}, which we call a {\it material G-structure} \cite{epselz}.

A principal-bundle connection in a material G-structure is called a {\it material connection} \cite{noll, wang}. The Christoffel symbols of a material connection can be expressed in terms of the matrices conveying the material isomorphims as
\begin{equation}
\Gamma^I_{JK}= P^I\;P^{-\alpha}_{\;\;\;J,K}.
\end{equation}
If the structure (symmetry) group is discrete, the material connection is unique. For the case of a continuous group, the material connection is not unique. The problem of {\it local material homogeneity} can be boiled down to answering the question as to whether a torsion-free material connection exists within the given material G-structure. Hence the importance of material connections in applications to theories of continuous distributions of defects in materials. According to the treatment of section \ref{sec:algebroidconnection} we should find all the information needed about material connections in the Lie algebroid associated with the material Lie groupoid.

\section{The material Lie algebroid}

For the sake of gathering a sort of physical picture, let us consider the case of a uniform body with a trivial symmetry group. Under this condition, the pointwise Lie algebra vanishes identically so that we can focus our attention on the material algebroid proper. Clearly, in the arrow picture, we have a single arrow (material isomorphism) for every ordered pair $(X,Y)$ of material points. Let $P^I_J(X,Y)$ denote the matrix components of this material isomorphism in a global Cartesian coordinate system in some reference configuration of the body. Adopting an archetypal point, this field of matrices can be expressed by the composition
\begin{equation}
P^I_J(X,Y)=P^I_\alpha(Y)\;P^{-\alpha}_{\;\;\;J}(X).
\end{equation}
Since we want to explore the vicinity of the identities, we leave the first argument $X$ (the tail of the arrows) fixed while we take derivatives with respect to the second argument $Y$ (the tip of the arrows) and we evaluate the derivatives as $Y \rightarrow X$. In this way, we obtain at each point $X$ the partial derivatives
\begin{equation}
\Gamma^I_{JK}= P^I_{\alpha,K}\;P^{-\alpha}_{\;\;\;J},
\end{equation}
which can be recognized as (minus) the Christoffel symbols of the (unique, in this case) material connection. Since this procedure mimics what one does to obtain a basis of the Lie algebra associated with a Lie group, we can conjecture intuitively that: {\it The material connections of a uniform body (together with a basis for the local Lie algebras) constitute a sort of `natural basis' for the material Lie algebroid}. In this rather intuitive way we have come full circle to show that, just as the natural concept of material groupoid encompasses all possible material G-structures, so too the concept of material Lie algebroid encompasses and generalizes the concept of material connection. From this perspective, we may say that, just as the Lie algebra of the material symmetry group expresses infinitesimal material symmetries at a point, so too the Lie algebroid of the material Lie groupoid encapsulates the idea of infinitesimal symmetries between different points.

\section{Conclusion}

Without claiming to have introduced any novel mathematical or  mechanical ideas, our objective in this article has been to provide a conceptual framework for a rigorous undertaking of the study of the geometrical underpinnings and implications of the constructs arising from the specification of the constitutive response of a material body. As pointed out in the introduction, the concept of material groupoid as a set of arrows joining points that share a material property is utterly intuitive. So, too, is the concept of material Lie algebroid, once we posit the smooth dependence of these arrows with respect to position within the body. The technical details may look daunting, but the physical and geometrical picture are crystal clear and speak for themselves.

\bigskip

{\bf Acknowledgment}: The financial support of the Natural Sciences and Engineering Research Council of Canada and of the Consejo Superior de Investigaciones Cient\'ificas of Spain is gratefully acknowledged.

\end{document}